\documentclass{amsart}

\usepackage{indentfirst}
\usepackage{amsmath}
\usepackage{amssymb}
\usepackage[dvipdfmx]{graphicx}
\usepackage{color}
\usepackage{amsthm}
\usepackage{mathrsfs}
\usepackage{pinlabel}
\usepackage{float}
\usepackage{amsopn}

\theoremstyle{plain}
 \newtheorem{theo}{Theorem}[section]

\theoremstyle{definition}

\theoremstyle{remark}

\newcommand{\liesl}{\mathfrak{sl}(2,\mathbb{C})}

\title{Periodicity Properties of the Colored Jones Polynomial}
\author{Hiroki Murakami}
\date{}

\begin{document}

\maketitle

\begin{abstract}
The "color" in the colored Jones polynomial is an integer parameter.
 In this paper, a periodic pattern of the values of the colored Jones polynomial at the second and the third roots of unity is found. If we substitute $-1$ to the colored Jones polynomial, the value is alternately $1$ or the determinant of the given link. When it comes to substituting the primitive third root of unity, another periodic pattern appears.
\end{abstract}

\section{Introduction}

The colored Jones polynomial is a framed link invariant.
It is a sequence of Laurent polynomials in one variable with integer coefficients and indexed by a positive integer $N$, the color.
Let $J_{K,N}(q)$ be the $N-$colored Jones polynomial of a link $K$ with indeterminate $q$.
Note that $J_{K,2}(q)$ coincides with the ordinary Jones polynomial. 
The main theorems of this paper are the following.

\begin{theo}
 For any $0$-framed knot $K$ and natural number $N$, we have
  \[ |J_{K,N}(-1)| = 
  \left \{
   \begin{array}{cl}
    1 & (N \hspace{1mm} \mathrm{is} \hspace{1mm} odd) \\
    \det K & (N \hspace{1mm} \mathrm{is} \hspace{1mm} even).
   \end{array}
  \right.
 \]
 Here $\det K = J_{K,2}(-1)$ is the determinant of $K$.
\end{theo}

Note that $J_{K,2}(-1) = \Delta _K(-1)$, where $\Delta_K(t)$ is the Alexander polynomial of $K$.
Along these lines, we also have

\begin{theo}
For $l \in \mathbb{Z}$, we obtain
 \[
 J_{K,N}\left( e^{\frac{2\pi \sqrt{-1}}{3}} \right) = \left \{
  \begin{array}{rl}
   0 & (N=6l) \\
   1 & (N=6l+1) \\
   1 & (N=6l+2) \\
   0 & (N=6l+3) \\
   -1 & (N=6l+4) \\
   -1 & (N=6l+5).
  \end{array}
 \right.
 \]
\end{theo}

\noindent
{\bf Organization.}
In section $2$, we describe the so-called cabling formula of the colored Jones polynomial. In section $3$, we prove our main theorems.

\noindent
{\bf Acknowledgements.}
The author should like to express his gratitude to Tam\'as K\'alm\'an for his constant encouragement and many pieces of helpful advice. 
He should also like to thank Keiji Tagami for his help proving Theorem1.1 and kind comments on results in this paper.


\section{Colored Jones polynomial}

Let $J_{K,N}(q)$ be the $N$-colored Jones polynomial with indeterminate $q$. We state a useful formula first to prove our results.

\begin{theo}[\cite{KM}, Theorem 4.15] \label{cabling}
 Let $K$ be a framed knot. Then 
 \[ J_{K, N}(q) = \sum_{j=0}^{ \lfloor \frac{N-1}{2} \rfloor} (-1)^j \binom{N-1-j}{j} J_{K^{N-1-2j}}(q), \]
 where $K^{N-1-2j}$ denotes the $0$-framed $N-1-2j$ cable of $K$. This formula is called the cabling formula.
\end{theo}

This formula says that the $N$-colored Jones polynomial can be computed using the ordinary Jones polynomial of cables.


\section{Value at Some Roots of Unity}

In this section, we examine the evaluation of the colored Jones polynomial at the second and the third roots of unity.
Note that our results are only proven in the case of knots, that is,  1-component links.

\subsection{Second Root of Unity}

Substituting the primitive second root of unity $-1$ to the colored Jones polynomial, we obtain the following. 

\begin{theo}
 Let $K$ be a $0$-framed knot. Then
 \[ |J_{K,N}(-1)| = 
  \left \{
   \begin{array}{cl}
    1 & (N \hspace{1mm} \mathrm{is} \hspace{1mm} odd) \\
    \det K & (N \hspace{1mm} \mathrm{is} \hspace{1mm} even).
   \end{array}
  \right.
 \]
 
 \begin{proof}
 For all $N \in \mathbb{N}_{\geq 2}$, we have 
  \[ J_{K,N}(q) = \sum_{j=0}^{ \lfloor \frac{N-1}{2} \rfloor} (-1)^j \binom{N-1-j}{j} J_{K^{N-1-2j}}(q) \]
  by the formula in Theorem \ref{cabling}. Using some well-known properties of the Jones polynomial and the Alexander polynomial, we obtain
   
   \begin{eqnarray*}
    J_{K,N}(-1) &=& \sum_{j=0}^{\lfloor \frac{N-1}{2} \rfloor} (-1)^j \binom{N-1-j}{j} J_{K^{N-1-2j}}(-1) \\
    			&=& \sum_{j=0}^{\lfloor \frac{N-1}{2} \rfloor} (-1)^j \binom{N-1-j}{j} \Delta_{K^{N-1-2j}}(-1) \\
			&=& \sum_{j=0}^{\lfloor \frac{N-1}{2} \rfloor} (-1)^j \binom{N-1-j}{j} \Delta_K\Bigl( (-1)^{N-1-2j} \Bigr) \cdot \Delta_{T_{N-1-2j,0}} (-1),
   \end{eqnarray*}
   where $T_{p,q}$ denotes the $(p,q)$-torus knot. If $N-1-2j$ is neither $0$ nor $1$, then $\Delta_{T_{N-1-2j,0}}(-1) =0$ because $T_{N-1-2j,0}$ is split.
   Note that we use the satellite formula to compute the Alexander polynomial of the cable.
   
   \begin{description}
    \item[$(1) N-1-2j=0$] For $N-1-2j$ to be $0$, $N$ must be an odd number. Since the upper number of the binomial coefficient must exceed the lower number, we obtain an equality $0=n-1-j \geq j$.
    					   It is equivalent to $j=0$, thus
    					   \[ | J_{K,N}(-1) | = \left| \sum_{j=0}^{\lfloor \frac{N-1}{2} \rfloor} (-1)^j \binom{N-1-j}{j} \right|= \left| \sum_{j=0}^{\lfloor \frac{N-1}{2} \rfloor} (-1)^j \right| = 1. \]
    
    \item[$(2) N-1-2j=1$] For $N-1-2j$ to be $1$, $N$ must be an even number. Then 
    					   \begin{eqnarray*}
    					    | J_{K,N}(-1) | &=& \left| \sum_{j=0}^{\lfloor \frac{N-1}{2} \rfloor} (-1)^j \binom{N-1-j}{j} | \Delta_K(-1)| \right| \\
					                         &=&|\Delta_K(-1)| \left| \sum_{j=0}^{\lfloor \frac{N-1}{2} \rfloor} (-1)^j (j+1) \right| \\
					                         &=& |\Delta_K(-1)|
					   \end{eqnarray*}
					   because we have $j=0$ and $j=1$ in this case. 
   \end{description}
   
  \end{proof} 
\end{theo}

\subsection{Third Root of Unity}

We examine the values of $J_{K,N}(q)$ at another fixed value of $q$. The main result of this subsection is the following.

\begin{theo}
 Let $K$ be a knot. Then 
 \[
  J_{K,N}\Bigl( e^{\frac{2\pi \sqrt{-1}}{3}} \Bigr) = 
   \left \{
    \begin{array}{rl}
      0 & (N=6l), \\
   1 & (N=6l+1), \\
   1 & (N=6l+2), \\
   0 & (N=6l+3), \\
   -1 & (N=6l+4), \\
   -1 & (N=6l+5).
    \end{array}
   \right.
 \]
 
 \begin{proof}
  From the fact that $J_{K,2}\Bigl( e^{\frac{2\pi \sqrt{-1}}{3}} \Bigr) = 1$ (See \cite{L}, for example) and by the formula in Theorem \ref{cabling}, we have
  \begin{eqnarray*}
   J_{K,N}\Bigl( e^{\frac{2\pi \sqrt{-1}}{3}}\Bigr ) &=& \sum_{j=0}^{\lfloor \frac{N-1}{2} \rfloor} (-1)^j \binom{N-1-j}{j} J_{K^{N-1-2j},2} \Bigl( e^{\frac{2\pi \sqrt{-1}}{3}} \Bigr) \\
   								&=& \sum_{j=0}^{\lfloor \frac{N-1}{2} \rfloor} (-1)^j \binom{N-1-j}{j}.
  \end{eqnarray*}
    Thus it suffices to show the combinatorial identity
  
  	\[ \sum^{\lfloor \frac{n}{2} \rfloor}_{j=0} (-1)^j \binom{n-j}{j} = \left \{ \begin{array}{rl} 
   		     1 & (n=6l), \\
   		     1 & (n=6l+1), \\
   		     0 & (n=6l+2), \\
   		     -1 & (n=6l+3), \\
   		     -1 & (n=6l+4), \\
   		     0 & (n=6l+5).
   		    \end{array}
   	\right. \]
  Here we put $n=N-1$.
  Let the left hand side of the previous equation be $a_n$. To verify the statement, we construct a generating function for $a_n$.
  That is, we construct a formal power series $f(x)$ whose coefficient of $x^n$ is $a_n$. Since
  
  	\[
  	 \begin{array}{lcl}
  	  (1-x)^0 & = & \binom{0}{0} \vspace{2mm} \\
  	  (1-x)^1 & = & \binom{1}{0} - \binom{1}{1}x \vspace{2mm }\\
  	  (1-x)^2 & = & \binom{2}{0} -\binom{2}{1}x + \binom{2}{2}x^2 \vspace{2mm} \\
  	  (1-x)^3 & = & \binom{3}{0} -\binom{3}{1}x + \binom{3}{2}x^2 -\binom{3}{3}x^3 \vspace{2mm} \\
  	  	      & \vdots & \\
	  (1-x)^k &=& \binom{k}{0} - \binom{k}{1}x + \dots + (-1)^k \binom{k}{k}x^k \\
	  	& \vdots &
  	 \end{array}
  	\]
  and $\displaystyle a_n=\binom{n}{0}-\binom{n-1}{1}+\binom{n-2}{2}+\dots +(-1)^{\lfloor \frac{n}{2} \rfloor}
  \binom{n-\lfloor \frac{n}{2} \rfloor}{\lfloor \frac{n}{2} \rfloor}$, the coefficient of $x^n$ in
  \[ f(x) = (1-x)^0 +x(1-x)+x^2(1-x)^2 + \cdots = \sum^{\infty }_{n=0} x^n(1-x)^n \]
  is $a_n$. Since $f(x)$ can be expanded as
  	\begin{eqnarray*}
  	 f(x) &=& \frac{1}{1-x(1-x)} \\
  	 	&=& \frac{1+x}{1+x^3} \\
  	 	&=& (1+x) \sum^{\infty }_{k=0} (-1)^k x^{3k} \\
  	 	&=& \sum^{\infty }_{k=0} (-1)^k(x^{3k}+x^{3k+1}) \\
  	 	&=& \sum^{\infty }_{l=0} (x^{6l} + x^{6l+1} +0x^{6l+2} -x^{6l+3} -x^{6l+4} + 0x^{6l+5}),
  	\end{eqnarray*}
 we have the desired result.
 \end{proof}
\end{theo}



\end{document}